\documentclass[12pt,a4paper,twoside]{article}

\usepackage{amsmath,amsfonts,amsthm,amsopn,color,amssymb}
\usepackage{palatino}
\usepackage{graphicx}

\newcommand{\eps}{\varepsilon}

\newcommand{\R}{\mathbb{R}}

\newcommand{\RN}{{\mathbb{R}^N}}
\newcommand{\RT}{{\mathbb{R}^3}}

\renewcommand{\le}{\leqslant}
\renewcommand{\ge}{\geqslant}

\renewcommand{\d }{\delta }

\newcommand{\g }{\gamma }

\renewcommand{\l }{\lambda}
\newcommand{\n }{\nabla }

\renewcommand{\t}{\theta}

\newcommand{\G}{\Gamma}

\renewcommand{\H}{H^1(\RN)}
\newcommand{\Hr}{H^1_r(\RN)}
\renewcommand{\P}{{\cal P}}
\newcommand{\Ne}{\mathcal{N}}

\newcommand{\D }{{\mathcal D}^{1,2}(\RN)}
\newcommand{\irn }{\int_{\RN}}
\newcommand{\irt }{\int_{\RT}}

\def\bbm[#1]{\mbox{\boldmath $#1$}}

\newtheorem{theorem}{Theorem}[section]
\newtheorem{lemma}[theorem]{Lemma}

\newtheorem{remark}[theorem]{Remark}
\newtheorem{corollary}[theorem]{Corollary}
\renewenvironment{proof}{\noindent{\textbf{Proof\quad}}}{$\hfill\square$\vspace{0.2 cm}\\}
\newenvironment{proofmain}{\noindent{\textbf{Proof of Theorem  \ref{main}\quad}}}{$\hfill\square$\vspace{0.2 cm}\\}
\newenvironment{proofmain2}{\noindent{\textbf{Proof of Theorem \ref{main2}\quad}}}{$\hfill\square$\vspace{0.2 cm}\\}
\newenvironment{proofradial}{\noindent{\textbf{Proof of Theorem  \ref{th:radial}\quad}}}{$\hfill\square$\vspace{0.2 cm}\\}
\newenvironment{proofnon}{\noindent{\textbf{Proof of Theorem \ref{non}\quad}}}{$\hfill\square$\vspace{0.2 cm}\\}
\newenvironment{proofconon}{\noindent{\textbf{Proof of Corollary \ref{co:non}\quad}}}{$\hfill\square$\vspace{0.2 cm}\\}

\title{{\bf On the Schr\"odinger equation in $\RN$\\
under the effect of a general\\
nonlinear term}}

\author{A. Azzollini \thanks{Dipartimento di Matematica, Universit\`a degli
Studi di Bari,  Via E. Orabona 4, I-70125 Bari, Italy, e-mail:
{\tt azzollini@dm.uniba.it}}
 \; \& \;
A. Pomponio\thanks{Dipartimento di Matematica, Politecnico di
Bari, Via E. Orabona 4, I-70125 Bari, Italy, e-mail: {\tt
a.pomponio@poliba.it}}}

\date{}

\begin{document}

\maketitle

\begin{abstract}
In this paper we prove the existence of a positive solution to the
equation $-\Delta u + V(x)u=g(u)$ in $\RN,$ assuming the general
hypotheses on the nonlinearity introduced by Berestycki \&
Lions. Moreover we show that a minimizing problem, related to the existence of a ground state, has no solution.
\end{abstract}

\section{Introduction}

This paper deals with the following equation:
\begin{equation}\label{eq:V}
\left\{
\begin{array}{ll}
-\Delta u + V(x)u=g(u), & x\in\RN, \quad N\ge 3;
\\
u>0.
\end{array}
\right.
\end{equation}
An existence result of nontrivial solutions for this kind of
problem has been obtained by Rabinowitz \cite{R} assuming that $g$
is superlinear and subcritical at infinity and satisfies the
global growth Ambrosetti-Rabinowitz condition
    \begin{equation}\label{eq:ambr}
        \exists \mu>2 \hbox{ s.t. } 0<\mu \int_0^t g(s) \,d s\le   g(t)t,
        \hbox{ for all }t\in \R.
    \end{equation}
This condition is used to get the boundedness of Palais-Smale
sequences. In \cite{JT2}, L. Jeanjean and K. Tanaka have been able
to remove the hypothesis \eqref{eq:ambr} by means of an abstract tool,
consisting in a suitable approximating method (see
\cite[Theorem~1.1]{J}). However they preserved the condition on
the superlinear growth at infinity.

On the other hand, in the fundamental paper \cite{BL1}, Berestycki
\& Lions proved the existence of a ground state, namely a
solution which minimizes the action among all the solutions, for
the problem
\begin{equation}\label{eq}
\left\{
\begin{array}{l}
-\Delta u =g(u), \hbox{ in } \RN;
\\
u \in H^1(\RN), u\neq 0,
\end{array}
\right.
\end{equation}
under the following assumptions on the nonlinearity $g$:
\begin{itemize}
\item[({\bf g1})] $g\in C(\RN,\R)$, $g$ odd; \item[({\bf g2})]
$-\infty <\liminf_{s\to 0^+} g(s)/s\le \limsup_{s\to 0^+}
g(s)/s=-m<0$; \item[({\bf g3})] $-\infty <\limsup_{s\to +\infty}
g(s)/s^{2^*-1}\le 0$; \item[({\bf g4})] there exists $\zeta>0$
such that $G(\zeta):=\int_0^\zeta g(s)\,d s>0$;
\end{itemize}
here $2^*=2N/(N-2)$.

So it seems that the ``natural" assumptions on the nonlinearity do
not require the superlinearity at infinity.

In view of this result, the aim of this paper is to study the
problem \eqref{eq:V} preserving the same general assumptions of
\cite{BL1} on $g$. Moreover we assume the following hypotheses on
$V:$
\begin{itemize}
\item[({\bf V1})] $V\in C^1(\RN,\R)$ and $V(x)\ge 0$, for all
$x\in \RN$, and the inequality is strict somewhere; \item[({\bf
V2})] $\|(\n V(\cdot)\mid \cdot)^+\|_{N/2}<2S$; \item[({\bf V3})]
$\lim_{|x|\to \infty}V(x)=0$; \item[({\bf V4})] $V$ is radially
symmetric;
\end{itemize}
here $(\n V(x)\mid x)^+=\max \big\{(\n V(x)\mid x),0\big\}$ and
$S$ is the best Sobolev constant of the embedding
$\D\hookrightarrow L^{2^*}(\RN),$ namely
\[
S  = \inf_{u\in\mathcal{D}^{1,2}\setminus\{0\}}\frac{\|\n
u\|_2^2}{\|u\|_{2^*}^2}.
\]

Our first main result is the following
\begin{theorem}\label{main}
Assume that ({\bf g1-4}) and ({\bf V1-4}) hold, then the problem \eqref{eq:V}
possesses at least a radially symmetric solution.
\end{theorem}

Up to our knowledge, this is the first result on a problem as \eqref{eq:V}, where exactly the same general hypotheses of Berestycki \& Lions \cite{BL1} are assumed on the nonlinearity $g$.

As a consequence of Theorem \ref{main}, we can prove the following

\begin{theorem}\label{th:radial}
Assuming that ({\bf g1-4}) and ({\bf V1-4}) hold, then the problem
\eqref{eq:V} possesses a radial ground state solution, namely a
solution minimizing the action among all the nontrivial radial
solutions.
\end{theorem}

\begin{remark}
We can generalize the hypotheses ({\bf V1}) and ({\bf V3}), requiring:
\begin{itemize}
\item[({\bf V1})'] $V\in C^1(\RN,\R)$ and $V(x)\ge V_0>-m$, for all $x\in \RN$, and the inequality is strict somewhere;
\item[({\bf V3})'] $\lim_{|x|\to \infty}V(x)=V_0$;
\end{itemize}
and supposing that the function $\tilde g(s)=g(s)-V_0s$ satisfies ({\bf g1-4}).
\end{remark}

\begin{remark}
The geometrical hypotheses on the potential $V$ do not allow us to
use concentration-compactness arguments as in \cite{JT2}. As a
consequence, we have to require a symmetry property on $V$ to
prevent any possible loss of mass at infinity.
\end{remark}

In the second part of the paper, we are interested in solving a
minimization problem strictly related to the existence of a ground
state solution for \eqref{eq:V}. Let
    \begin{equation*}
        I(u)=\frac 12 \irn |\n u|^2 + V(x) u^2 -\irn G(u),\quad u\in\H.
    \end{equation*}
When the nonlinearity $g$ satisfies \eqref{eq:ambr}, a standard
method to look for the existence of a ground state solution for an
equation as \eqref{eq:V} is to study the minimizing problem
    \begin{equation}\label{eq:minimizing}
        I(\bar u)=\inf_{u\in\Ne} I(u),\quad \bar u \in \Ne,
    \end{equation}
where $\Ne$ is the Nehari manifold related to $I$. In \cite{R} the
geometrical assumption on $V$
    \begin{equation}\label{eq:geo}
        V(y)\le\lim_{x\to\infty} V(x),\quad\hbox{for all }
        y\in\R^N
    \end{equation}
is used to solve such a minimizing problem. Moreover, in
\cite{JT2} it has been proved that, assuming \eqref{eq:geo}, there
exists a ground state solution for \eqref{eq:V} also without the
condition \eqref{eq:ambr}. On the other hand, it is well known
(see for example \cite{BGM}) that if $g$ satisfies
\eqref{eq:ambr}, and \eqref{eq:geo} holds with the reverse
inequality, the minimizing problem \eqref{eq:minimizing} cannot
be solved. In fact, a contradiction argument deriving from the
comparison of the level $\eta:=\inf_{u\in\Ne} I(u)$ with $\eta_0
:= \inf_{u\in\Ne_{0}}I_{0}(u)$ (here $I_0$ and $\Ne_0$ are the
functional and the Nehari manifold of the problem at infinity) is
used and the Ambrosetti-Rabinowitz condition plays a fundamental
role. Actually, when we do not assume such a type of growth
condition on $g$, these arguments do not work any more. However, a
similar study can be done by replacing the Nehari manifold with a
more suitable one.\\
Indeed, if we define by
\begin{equation}\label{eq:P}
\P_0:=\left\{ u\in \H\setminus\{0\} \; \Big{|}\; \frac{N-2}2 \irn
|\n u|^2 = N \irn G(u)\right\},
\end{equation}
the Pohozaev manifold related to \eqref{eq}, and we set
$$
{\mathcal S_0}:=\{u\in\H\setminus\{0\}\mid u\hbox{ is a solution of \eqref{eq}}\},
$$
it is well known that ${\mathcal S_0}\subset \P_0.$ Moreover, in
\cite{S} it has been proved that $\P_0$ is a natural constraint
for the functional related to \eqref{eq}
\[
I_0(u)=\frac 12 \irn |\n u|^2 - \irn G(u),\quad u\in\H.
\]
In order to look for a
ground state solution, a very natural question arises: is the infimum $I_0|_{\P_0}$ achieved?

Following two different ways, Jeanjean \& Tanaka \cite{JT} and
Shatah \cite{S} have given a positive answer to this question,
showing that
\begin{equation}    \label{eq:SP}
\min_{u \in {\cal S}_0}I_0(u)=\min_{u \in {\cal P}_0}I_0(u).
\end{equation}
Inspired by these papers and observed that each solution of
\eqref{eq:V} satisfies the following Pohozaev identity:
\begin{equation}    \label{eq:poho}
\frac{N-2}{2}\irn |\n u|^2 + \frac N2 \irn V(x)u^2 +\frac 12 \irn
(\n V(x)\mid x)u^2 = N \irn G(u),
\end{equation}
we indicate with $\P$ the Pohozaev manifold related to
\eqref{eq:V}:
\begin{equation}\label{eq:PV}
\P=\left\{ u\in \H\setminus\{0\} \mid u \hbox{ satisfies
\eqref{eq:poho}} \right\}
\end{equation}
and we wonder if there exists a minimizer for $I|_\P$. 

Proceeding in analogy with \cite{BGM}, we get the following result

\begin{theorem}\label{non}
If we assume ({\bf g1-4}), ({\bf V1}), ({\bf V3}) and
\begin{itemize}
\item[({\bf V5})] $(\n V(x)\mid x)\le 0,$ for all $x\in\RN;$
\item[({\bf V6}] $NV(x)+(\n V(x)\mid x)\ge 0$ for all $x\in \RN$, and the inequality is strict somewhere,
\end{itemize}
then $b:=\inf_{u\in\P}I(u)$ is not a critical level for the
functional $I$.
\end{theorem}

Requiring something more, instead of ({\bf V6}), we have a further
information on the level $b$:
\begin{corollary}\label{co:non}
Assuming the same hypotheses of Theorem \ref{non}, if in ({\bf
V6}) we have the strict inequality, then the level $b$ is not
achieved as a minimum on the Pohozaev manifold $\P$.
\end{corollary}

{\bf Acknowledgement} The authors express their deep gratitude to
Prof. L.~Jeanjean for stimulating discussions and useful
suggestions.

\vspace{0.5cm}
\begin{center}
{\bf NOTATION}
\end{center}

\begin{itemize}
\item For any $1\le s< +\infty$, $L^s(\RT)$ is the usual Lebesgue space endowed with the norm
\[
\|u\|_s^s:=\irt |u|^s;
\]
\item $\H$ is the usual Sobolev space endowed with the norm
\[
\|u\|^2:=\irt |\n u|^2+ u^2;
\]
\item $\D$ is completion of $C_0^\infty(\RT)$ with respect to the norm
\[
\|u\|_{\D}^2:=\irt |\n u|^2;
\]
\item for any $r>0,$ $x\in\RT$ and $A\subset \RT$
\begin{align*}
B_r(x) &:=\{y\in\RT\mid |y-x|\le r\},\\
B_r    &:=\{y\in\RT\mid |y|\le r\},\\
A^c    &:= \RT\setminus A.
\end{align*}
\end{itemize}

\section{The existence result}

The aim of this section is to prove Theorems \ref{main} and \ref{th:radial}.

Set
    \begin{equation*}
        H^1_r(\R^N):=\{u\in\H\mid u \hbox{ is radial}\}
    \end{equation*}
and, following \cite{BL1}, define $s_0:=\min\{s\in
[\zeta,+\infty[\;\mid g(s)=0\}$ ($s_0=+\infty$ if $g(s)\neq 0$ for
any $s\ge\zeta$) and set $\tilde g:\R\to\R$ the function such that
    \begin{equation}\label{eq:tilde}
      \tilde g(s)=\left\{
        \begin{array}{ll}
                g(s) &\hbox{ on } [0,s_0];
                \\
                0 &\hbox{ on } \R_+\setminus [0,s_0];
                \\
                -\tilde g(-s) &\hbox{ on } \R_-.
      \end{array}
      \right.
    \end{equation}
By the strong maximum principle, a solution of \eqref{eq:V} with
$\tilde g$ in the place of $g$ is a solution of \eqref{eq:V}. So
we can suppose that $g$ is defined as in \eqref{eq:tilde}, so that
({\bf g1}), ({\bf g2}), ({\bf g4}) and then the following limit
    \begin{equation}\label{eq:limg}
        \lim_{s\to\infty} \frac{|g(s)|}{|s|^{2^*-1}}=0
    \end{equation}
hold. Moreover, we set for any $s\ge 0,$
    \begin{align*}
        g_1(s) & :=(g(s)+ms)^+,
        \\
        g_2(s) & :=g_1(s)-g(s),
    \end{align*}
and we extend them as odd functions.\\
Since
    \begin{align}
        \lim_{s\to 0}\frac{g_1(s)}{s} &= 0,\label{eq:lim1}\\
        \lim_{s\to\infty}\frac{g_1(s)}{|s|^{2^*-1}}&=0,\label{eq:lim2}
    \end{align}
and
    \begin{equation}
        g_2(s) \ge ms,\quad  \forall s\ge 0,\label{eq:g2}
    \end{equation}
by some computations, we have that for any $\eps>0$ there exists
$C_\eps>0$ such that
    \begin{equation}
        g_1(s) \le C_\eps s^{2^*-1}+\eps g_2(s),\quad  \forall
        s\ge0\label{eq:g1g2}.
    \end{equation}
If we set
    \begin{equation*}
        G_i(t):=\int^t_0g_i(s)\,ds,\quad i=1,2,
    \end{equation*}
then, by \eqref{eq:g2} and \eqref{eq:g1g2}, we have
    \begin{equation}
         G_2(s) \ge \frac m 2 s^2,\quad  \forall s\in\R\label{eq:G2}
    \end{equation}
and for any $\eps>0$ there exists $C_\eps>0$ such that
    \begin{equation}
        G_1(s) \le \frac {C_\eps} {2^*} |s|^{2^*}+\eps G_2(s),\quad  \forall
        s\in\R\label{eq:G1G2}.
    \end{equation}

Using an idea from \cite{J}, we look for bounded Palais-Smale
sequences of the following perturbed functionals
\begin{equation*}
I_\l(u)=\frac 12 \irn |\n u|^2 + V(x)u^2 + 
\irn G_2(u) - \l\irn G_1(u),
\end{equation*}
for almost all $\l$ near $1$. Then we will deduce the existence of
a non-trivial critical point $v_\l$ of the functional $I_\l$ at
the mountain pass level. Afterward, we study the convergence of
the sequence $(v_\l)_\l$, as $\l$ goes to 1 (observe that
$I_1=I$).

We will apply the following slight modified version of
\cite[Theorem~1.1]{J} (see \cite{J2}).

    \begin{theorem}\label{th:J}
        Let $\big(X,\|\cdot\|\big)$ be a Banach space and
        $J\subset\R^+$ an interval.
Consider the family of $C^1$ functionals on $X$
    \begin{equation*}
        I_\l(u)=A(u)- \l B(u),\quad\forall\l\in J,
    \end{equation*}
with $B$ nonnegative and either $A(u)\to + \infty$ or
$B(u)\to+\infty$ as $\|u\|\to\infty.$\\
For any $\l\in J$ we set
    \begin{equation}\label{eq:gamma}
        \Gamma_\l:=\{\gamma\in C([0,1],X)\mid \gamma(0)=0\neq\gamma(1),
        I_\l(\gamma(1))< 0,\}.
    \end{equation}
If for every $\l\in J$ the set $\G_\l$ is nonempty and
    \begin{equation}\label{eq:cl}
        c_\l:=\inf_{\gamma\in\Gamma_\l}\max_{t\in[0,1]}
        I_\l(\gamma(t)) > I_\l(v),
    \end{equation}
then for almost every $\l\in J$ there is a sequence $(v_n)_n
\subset X$ such that
    \begin{itemize}
        \item[(i)] $(v_n)_n$ is bounded;
        \item[(ii)] $I_\l(v_n)\to c_\l$;
        \item[(iii)] $(I_\l)'(v_n)\to 0$ in the dual $X^{-1}$ of $X$.
    \end{itemize}
    \end{theorem}
In our case, $X=\Hr$ and
    \begin{align*}
        A(u) &:=\frac 12 \irn |\n u|^2 + V(x)u^2 
        + \irn G_2(u),\\
        B(u) &:=\irn G_1(u).
    \end{align*}

In order to apply Theorem \ref{th:J}, we have just to define a
suitable interval $J$ such that $\Gamma_\l\neq\emptyset$, for any $\l\in J$, 
and \eqref{eq:cl} holds.\\
Observe that, according to \cite{BL1}, there exists a function
$z\in\Hr$ such that
    $$\irn G_1(z)-\irn G_2(z)=\irn G(z)>0.$$
Then there exists $0<\bar\d<1$ such that
    \begin{equation}\label{eq:d}
        \bar\d\irn G_1(z)-\irn G_2(z)>0.
    \end{equation}
We define $J$ as the interval $[\bar \delta , 1].$

\begin{lemma}\label{le:Gamma}
$\Gamma_\l\neq\emptyset$, for any $\l\in J$.
\end{lemma}

\begin{proof}
Let $\l\in J$. Set $\bar\t>0$ sufficiently large and $\bar
z=z(\cdot /
\bar\t)$.\\
Define $\g:[0,1]\to\Hr$ in the following way
    \begin{equation*}
        \g(t)=\left\{\begin{array}{ll}
            0,& \hbox{if }t=0,
            \\
            \bar z^t=\bar z(\cdot/t), & \hbox{if }0<t\le1.
        \end{array}
        \right.
    \end{equation*}
It is easy to see that $\g$ is a continuous path from $0$ to $\bar
z.$ Moreover, we have that
    \begin{multline*}
        I_\l(\g(1)) \le \frac {\bar\t^{N-2}} 2 \irn |\n z|^2 
        +\frac{\bar\t^N} 2 \irn V\left(\bar\t x\right) z^2 
        \\
        +\bar\t^N
        \left(\irn G_2(z) - \bar\d\irn G_1(z)\right)
    \end{multline*}
and then, by \eqref{eq:d}, ({\bf V3}) and the Lebesgue theorem,
for a suitable choice of $\bar\t$, certainly $\g\in\G_\l.$
\end{proof}

\begin{lemma}\label{le:cl}
$c_\l>0$ for all $\l\in J.$
\end{lemma}

\begin{proof}
Observe that for any $u\in\Hr$ and $\l\in J$, using \eqref{eq:G2}
and \eqref{eq:G1G2} for $\eps < 1$, we have
    \begin{align*}
        I_\l(u) & \ge \frac 1 2 \irn |\n u|^2 + V(x)u^2
        +\irn G_2 (u) -\irn G_1 (u)\\
                  & \ge \frac 1 2 \irn |\n u|^2  + (1-\eps) \frac
                  m 2 \irn u^2 - \frac{C_\eps}{2^*} \irn
                  |u|^{2^*}
    \end{align*}
and then, by Sobolev embeddings, we conclude that there exists
$\rho >0$ such that for any $\l\in J$ and $u\in\Hr$ with $u\neq 0$ and
$\|u\|\le\rho,$ it results $I_\l(u)>0.$ In particular, for any
$\|u\|=\rho,$ we have $I_\l(u) \ge \tilde c >0.$ Now fix $\l\in J$
and $\g\in\G_\l.$ Since $\g(0)=0$ and $I_\l(\g(1))< 0$,
certainly $\|\g(1)\|
> \rho.$ By continuity, we deduce that there exists $t_\g\in
]0,1[$ such that $\|\g(t_\g)\|=\rho.$ Therefore, for any $\l\in
J,$
\begin{equation}\label{eq:cll}
c_\l\ge \inf_{\g\in\G_\l} I_\l(\g(t_\g)) \ge \tilde c >0.
\end{equation}
\end{proof}
We present a variant of the Strauss' compactness lemma \cite{Str}
(see also \cite[Theorem A.1]{BL1}), whose proof is similar to that
contained in \cite{BL1}. It will be a fundamental tool in our
arguments:
\begin{lemma}\label{le:str}
Let $P$ and $Q:\R\to\R$ be two continuous functions satisfying
\begin{equation*}
\lim_{s\to\infty}\frac{P(s)}{Q(s)}=0,
\end{equation*}
$(v_n)_n,$ $v$ and $z$ be measurable functions from $\RN$ to $\R$
such that
\begin{align*}
&\sup_n\irn | Q(v_n(x))z|\,dx <+\infty,
\\ 
&P(v_n(x))\to v(x) \:\hbox{a.e. in }\RN. 
\end{align*}
Then $\|(P(v_n)-v)z\|_{L^1(B)}\to 0$, for any bounded Borel set
$B$.

Moreover, if we have also
\begin{align*}
\lim_{s\to 0}\frac{P(s)}{Q(s)} &=0,\\ 
\lim_{x\to\infty}\sup_n |v_n(x)| &= 0, 
\end{align*}
then $\|(P(v_n)-v)z\|_{L^1(\RN)}\to 0.$
\end{lemma}

In analogy with the well-known compactness result in \cite{BL2},
we state the following result

\begin{lemma}\label{le:PS}
For any $\l \in J$, each bounded Palais-Smale sequence for the functional $I_\l$ admits a convergent subsequence.
\end{lemma}

\begin{proof}
Let $\l\in J$ and $(u_n)_n$ be a bounded (PS) sequence for $I_\l$,
namely
    \begin{align}\label{eq:PS}
        &(I_\l(u_n))_n \hbox{ is bounded },\nonumber\\
        &\lim_n (I_\l)'(u_n)= 0 \hbox{ in } (\Hr)'.
    \end{align}
Up to a subsequence, we can suppose that there exists $u\in\Hr$
such that
    \begin{equation}\label{eq:weak}
        u_n\rightharpoonup u\;\hbox{weakly in }\Hr
    \end{equation}
and
    \begin{equation}\label{eq:aeconv}
        u_n(x)\to u(x)\;\hbox{a.e. in }\RN.
    \end{equation}
By weak lower semicontinuity we have:
\begin{align}
\irn |\n u|^2 \le &\liminf_n \irn |\n u_n|^2; \label{eq:semi1}
\\
\irn V(x) u^2 \le &\liminf_n \irn V(x)u_n^2. \label{eq:semi2}
\end{align}
If we apply Lemma \ref{le:str} for $P(s)=g_i(s)$, $i=1,2,$ $Q(s)=
|s|^{2^*-1},$ $(v_n)_n=(u_n)_n,$ $v=g_i(u),$ $i=1,2$ and $z\in
C^\infty_0(\RN),$ by \eqref{eq:limg}, \eqref{eq:lim2} and
\eqref{eq:aeconv} we deduce that
    \begin{equation*}
        \irn g_i(u_n)z\to\irn g_i(u)z\quad i=1,2.
    \end{equation*}
As a consequence, by \eqref{eq:PS} and \eqref{eq:weak} we deduce
$(I_\l)'(u)=0$ and hence
\begin{equation}    \label{eq:solo}
\irn |\n u|^2+ V(x)u^2= \irn \l g_1(u)u - g_2(u)u.
\end{equation}
\\
If we apply Lemma \ref{le:str} for $P(s)=g_1(s)s,$ $Q(s)= s^2+
|s|^{2^*},$ $(v_n)_n=(u_n)_n,$ $v=g_1(u)u,$ and $z=1,$ by
\eqref{eq:limg}, \eqref{eq:lim2}, \eqref{eq:aeconv} and the well
known Strauss' radial lemma (see \cite{Str}) we deduce that
    \begin{align}
        \irn g_1(u_n)u_n \to \irn g_1(u)u.\label{eq:convg}
    \end{align}
Moreover, by \eqref{eq:aeconv} and Fatou's lemma
    \begin{align}
        \irn g_2(u)u\le &\liminf_n \irn g_2(u_n)u_n.\label{eq:convg2}
    \end{align}
By \eqref{eq:solo}, \eqref{eq:convg} and \eqref{eq:convg2}, and
since $\langle (I_\l)'(u_n),u_n\rangle\to0$
\begin{align}
\limsup_n \irn |\n u_n|^2 + V(x)u_n^2 &=\limsup_n \left[\l \irn
g_1(u_n)u_n -\irn g_2(u_n)u_n\right] \nonumber
\\
&\le\l \irn g_1(u)u - \irn g_2(u)u \nonumber\\
&=\irn |\n u|^2 + V(x)u^2. \label{eq:limsup}
\end{align}
By \eqref{eq:semi1}, \eqref{eq:semi2} and \eqref{eq:limsup}, we get
\begin{align}
\lim_n \irn |\n u_n|^2= &\irn |\n u|^2 \label{eq:sf1},
\\
\lim_n \irn V(x) u_n^2= &\irn V(x) u^2, \nonumber 
\end{align}
hence
\begin{equation}    \label{eq:sfg2}
\lim_n \irn g_2(u_n)u_n= \irn g_2(u)u.
\end{equation}
Since $g_2(s)s=ms^2 +q(s)$, with $q$ a positive and continuous function, by Fatou's Lemma we have
\begin{align*}
\irn q(u) \le &\liminf_n \irn q(u_n);
\\
\irn u^2 \le &\liminf_n \irn u_n^2.
\end{align*}
These last two inequalities and \eqref{eq:sfg2} imply that, up to
a subsequence,
\[
\irn  u^2 =\lim_n \irn u_n^2,
\]
which, together with \eqref{eq:sf1}, shows that $u_n \to u$ strongly in $\Hr$.
\end{proof}

\begin{lemma}\label{le:ul}
For almost every $\l\in J$, there exists $u^\l\in \Hr$, $u^\l\neq
0$, such that $(I_\l)'(u^\l)=0$ and $I_\l(u^\l)= c_\l$.
\end{lemma}

\begin{proof}
By Theorem \ref{th:J}, for almost every $\l\in J$, there exists a
bounded sequence $(u^\l_n)_n\subset\Hr$ such that
    \begin{align}
        I_\l(u^\l_n)&\to c_\l;\label{eq:conv}\\
        (I_\l)'(u^\l_n)&\to 0\;\hbox{in } (\Hr)'.\label{eq:p-s}
    \end{align}
Up to a subsequence, by Lemma \ref{le:PS}, we can suppose that
there exists $u^\l\in\Hr$ such that $u^\l_n \to u^\l$ in $\Hr$. By
Lemma \ref{le:cl}, \eqref{eq:conv} and \eqref{eq:p-s} we conclude.
\end{proof}

Now we are able to provide the proof of our main result:

\begin{proofmain}
By Lemma \ref{le:ul}, we are allowed to consider a suitable
$\l_n\nearrow 1$ such that for any $n\ge 1$ there exists
$v_n\in\Hr\setminus\{0\}$ satisfying
\begin{align}
I_{\l_n}(v_n)&= c_{\l_n},\label{eq:mp}
\\
(I_{\l_n})'(v_n)&=0\;\hbox{in }(\Hr)'.      \label{eq:sol}
\end{align}
We want to prove that $(v_n)_n$ is a bounded Palais-Smale sequence for $I$ at the level $c:=c_1$.\\
By standard argument, since $\Hr$ is a natural constraint, we have
that $v_n$ is a weak solution of the problem
    \begin{equation*}
        -\Delta w + V(x) w + g_2(w) - \l_n g_1(w) = 0
    \end{equation*}
and it satisfies the Pohozaev equality
    \begin{multline}\label{eq:Poho}
        \irn |\n v_n|^2 + \frac N {N-2} \irn V(x) v_n^2
        + \frac 1 {N-2}\irn (\n V(x)\mid x) v_n^2
        \\
        +\frac {2N} {N-2}
        \irn G_2(v_n)-\l_n G_1(v_n)=0.
    \end{multline}
Therefore, by \eqref{eq:mp}, \eqref{eq:sol} and \eqref{eq:Poho} we
have that the following system holds
    \begin{equation}
        \left\{
            \begin{array}{l}
                \frac 1 2 (\alpha_n + \beta_n) + \gamma_{2,n} - \l_n
                \gamma_{1,n} = c_{\l_n},\label{eq:sys1}
                \\
                \alpha_n + \beta_n + \d_{2,n} - \l_n \d_{1,n} = 0,
                \\
                \alpha_n + \frac N {N-2} \beta_n + \frac 1 {N-2}\eta_n +\frac {2N} {N-2}
                \gamma_{2,n} - \frac {2N} {N-2}\l_n \gamma_{1,n}=0,
            \end{array}
        \right.
    \end{equation}
where
    \begin{eqnarray*}
        \alpha_n = \irn |\n v_n|^2, & \displaystyle\beta_n=\irn V(x) v_n^2,& \eta_n=\irn (\n
        V(x)\mid x) v_n^2,
        \\
        \gamma_{i,n}=\irn G_i(v_n),   & \displaystyle\d_{i,n} =\irn g_i(v_n)v_n, &
        i=1,2.
    \end{eqnarray*}
By the first and the third of the system we get
    \begin{equation*}
        \frac 1 N \alpha_n - \frac 1 {2N} \eta_n = c_{\l_n}
    \end{equation*}
and then, using Holder inequality, by $({\bf V2})$ and the boundedness of $(c_{\l_n})_n$
(indeed the map $\l \mapsto c_\l$ is non-increasing), we have
    \begin{equation}\label{eq:bound}
        \alpha_n\le C \hbox{ for all }n\ge 1.
    \end{equation}
By the second of the system we have
    \begin{equation*}
        \d_{2,n} - \l_n \d_{1,n} = - \alpha_n  - \beta_n\le 0
    \end{equation*}
and then by \eqref{eq:g1g2}, there
exists $0<\eps<1$ and $C_\eps>0$ such that
    \begin{equation*}
        \d_{2,n} \le \d_{1,n} \le C_\eps \irn |v_n|^{2^*}+\eps
        \d_{2,n}.
    \end{equation*}
Therefore, by the Sobolev embedding $\D\hookrightarrow
L^{2^*}(\RN)$
    \begin{equation*}
        (1-\eps) \d_{2,n} \le C_\eps \irn |v_n|^{2^*}\le C
        \alpha_n^{2^*}
    \end{equation*}
and then, by \eqref{eq:bound}, $\d_{2,n}$ is bounded. By \eqref{eq:g2} and
\eqref{eq:bound} we deduce that $(v_n)_n$ is bounded in $\Hr$.\\
Up to a
subsequence, there exists $v\in\Hr$ such that
    \begin{equation}\label{eq:weak2}
        v_n\rightharpoonup v\;\hbox{weakly in }\Hr.
    \end{equation}
By \eqref{eq:sol}, we have that
    \begin{equation*}
        I'(v_n)=(I_{\l_n})'(v_n) + (\l_n -1) g_1(v_n)= (\l_n-1)
        g_1(v_n)
    \end{equation*}
so 
\begin{equation}	\label{eq:nome}
I'(v_n)\to 0 \quad \hbox{in }(\Hr)',
\end{equation}
if $\left(g_1(v_n)\right)_n$ is bounded in $\big(\Hr\big)'.$
But this is true by the Banach-Steinhaus theorem, since Lemma
\ref{le:str} implies that for any $z\in\Hr$
\begin{equation*}
\irn g_1(v_n)z \to \irn g_1(v)z.
\end{equation*}
Moreover, from \eqref{eq:mp} and the boundedness of $(v_n)_n$, we deduce that
    \begin{equation}\label{eq:mou}
        I(v_n)=I_{\l_n}(v_n) + (\l_n -1) \irn G_1(v_n)\to c.
    \end{equation}
Therefore, by \eqref{eq:nome} and \eqref{eq:mou}, $(v_n)_n$ is a Palais-Smale sequence for the functional $I$ and so, by Lemma \ref{le:PS}, $v$ is a nontrivial mountain pass type solution for \eqref{eq:V}.
\\
To conclude, observe that, by standard arguments, we can use the
strong maximum principle to get $v>0.$

\end{proofmain}
In order to prove Theorem \ref{th:radial} we set
    \begin{align*}
        {\mathcal S_r}&:=\left\{u\in\Hr\setminus\{0\}\mid
        I'(u)=0\right\},\\
        \sigma_r&:=\inf_{u\in {\mathcal S_r}} I(u).
    \end{align*}
\begin{lemma}\label{le:s}
    We have $\sigma_r>0.$
\end{lemma}
\begin{proof}
    By \eqref{eq:g1g2} and Sobolev embedding we have that for any $u\in{\mathcal S_r}$
        \begin{align*}
            \|\n u\|_2^2&\le \irn |\n u|^2 + V(x) u^2 +
            (1-\eps)\irn g_2(u)u\\
            &\le C_\eps \|u\|_{2^*}^{2^*}\le C \|\n u\|_2^{2^*}
        \end{align*}
    where $0<\eps<1$ and $C_\eps,C>0.$
    So we deduce that
        \begin{equation*}
            \inf_{u\in{\mathcal S_r}} \|\n u\|_2>0.
        \end{equation*}
    Now, since ${\mathcal S_r}\subset \P$ (see \eqref{eq:PV}), for any $u\in{\mathcal S_r}$ by ({\bf V2}) we have
        \begin{equation*}
            I(u)= \frac 1N\irn |\n u|^2 -\frac 1{2N} \irn (\n V(x)\mid x) u^2\ge C >0 .
        \end{equation*}
\end{proof}
Finally we provide the following

\begin{proofradial}
Let $(u_n)_n\in{\mathcal S_r}$ such that $I(u_n)\to \sigma_r.$
Arguing as in the proof of Theorem \ref{main} we have that the
sequence is bounded. By Lemma \ref{le:PS}, there exists $u\in\Hr$
such that $u_n\to u$ in $\Hr$ and then the conclusion follows by Lemma \ref{le:s}.
\end{proofradial}

\section{The nonexistence result}

In this section we give the proof of Theorem \ref{non} and we will assume that $V$ satisfies hypotheses ({\bf V1}), ({\bf V3}), ({\bf V5-6}).

Let us show that the functional $I$ is bounded below on the
manifold $\P$:
\begin{lemma}\label{le:PV}
For all $u \in \P$, $I(u)>0$.
\end{lemma}

\begin{proof}
It is easy to see that for any $u\in \P$, by ({\bf V2}) we get
\begin{equation}    \label{eq:iv}
I(u)= \frac 1N\irn |\n u|^2 -\frac 1{2N} \irn (\n V(x)\mid x) u^2>
0.
\end{equation}
\end{proof}

By Lemma \ref{le:PV}, we can define
\begin{equation}    \label{eq:mv}
b=\inf_{u \in \P}I(u)\ge 0.
\end{equation}

\begin{lemma}\label{le:p-pv}
Let $w\in \H$ be such that $\irn G(w)>0$. Then there exists
$\bar\t>0$ such that $w^{\bar\t}=w(\cdot/\bar\t)\in \P$. In
particular this result is true for any $w\in \P_0$ (see
\eqref{eq:P}).
\end{lemma}

\begin{proof}
For any $\t>0$, we set
\[
f(\t):=I(w^\t) =\frac{\t^{N-2}}{2} \irn |\n w|^2 +\frac{\t^N}{2}
\irn V(\t x) w^2 - \t^N \irn G(w).
\]
By the Lebesgue theorem and ({\bf V3}), we get
\[
\lim_{\t\to +\infty}\irn V(\t x) w^2=0,
\]
and so
\[
\lim_{\t\to +\infty}f(\t)=-\infty.
\]
We argue that there exists $\bar\t>0$ such that $f'(\bar\t)=0$:
hence $w^{\bar\t}\in \P$.
\end{proof}

Let $w\in \P_0$. For any $y\in \RN$, we set $w_y:=w(\cdot
-y)\in{\mathcal P_0}$. Set $\t_y>0$ such that $\tilde
w_y=w_y(\cdot/\t_y)\in \P$.

\begin{lemma}\label{le:t}
We have $\lim_{|y|\to \infty}\t_y=1$.
\end{lemma}

\begin{proof}
{\sc Step 1}: $\limsup_{|y|\to\infty}\t_y<+\infty.$
\\
Suppose, by contradiction, that $\t_{y_n} \to +\infty$, for
$|y_n|\to \infty$.
\\
For any $y\in \RN$, we have
\begin{equation}\label{eq:iwt}
I(\tilde w_y) =\frac{\t_y^{N-2}}{2} \irn |\n w|^2
+\frac{\t^N_y}{2} \irn V(\t_y x) \,w^2\!\left(\!x-\frac{y}{\t_y}
\!\right) - \t^N_y \irn G(w).
\end{equation}
Let us show that
\begin{equation}    \label{eq:wt}
\lim_{|y|\to \infty}\irn V(\t_y x) \,w^2\!\left(\!x-\frac{y}{\t_y}
\!\right)=0.
\end{equation}
Indeed we have
\begin{align}
\irn V(\t_y x) \,w^2\!\left(\!x-\frac{y}{\t_y} \!\right)
&=\int_{B_r} V(\t_y x) \,w^2\!\left(\!x-\frac{y}{\t_y} \!\right)
+\int_{B_r^c} V(\t_y x) \,w^2\!\left(\!x-\frac{y}{\t_y} \!\right) \nonumber
\\
&\le \sup_{x\in\RN}V(x)\int_{B_r(-y/\t_y)}w^2 
+\sup_{x\in B_r^c}V(\t_y x)\,\|w\|_2^2. \label{eq:dai}
\end{align}
By the absolute continuity of the Lebesgue integral, for any
$\eps>0$ there exists $\bar r>0$ such that, for any $r<\bar r$ and
for any $y\in \RN$, we get
\begin{equation}	\label{eq:ihih}
\sup_{x\in\RN}V(x) \cdot \int_{B_r(-y/\t_y)}w^2 \le \eps.
\end{equation}
Therefore, since we are supposing that $\t_{y_n} \to +\infty$, as
$|y_n|\to \infty$, by \eqref{eq:dai}, \eqref{eq:ihih} and ({\bf V2}), we get \eqref{eq:wt}. As a consequence, by
\eqref{eq:iwt} and \eqref{eq:wt}, we infer that $I(\tilde
w_{y_n})\to -\infty$, as $|y_n|\to \infty$, and we get a
contradiction with Lemma \ref{le:PV}.
\\
\
\\
{\sc Step 2}: $\lim_{|y|\to \infty}\t_y=1$.
\\
Since $w\in \P_0$ and $\tilde w_y\in \P$, we get
\begin{multline}\label{eq:ty}
N(\t_y^2-1)\irn \! G(w)\\
=\frac{\t_y^2}{2} \irn \big[N V(\t_y
x+y) +(\n V(\t_y x+y)\mid(\t_y x+y))\big]w^2.
\end{multline}
By ({\bf V3}), ({\bf V5}) and ({\bf V6}), using the dominated
convergence and the conclusion of the Step 1, the right hand side
in \eqref{eq:ty} goes to zero as $|y|\to \infty$, and so the lemma
is proved.
\end{proof}

We set (see \eqref{eq:SP})
\[
b_0:=\min_{u\in{\cal S}_0}I_0(u)=\min_{u\in \P_0}I_0(u).
\]

\begin{lemma}\label{le:<=}
$b\le b_0$.
\end{lemma}

\begin{proof}
Let $w\in \H$ be a ground state solution of \eqref{eq}. Then $w\in
\P_0$ and $I_0(w)=b_0$. For any $y\in \RN$, we set $w_y=w(\cdot
-y)$. By the invariance by translations of \eqref{eq}, we have
that $w_y\in \P_0$ and $I_0(w_y)=b_0$. By Lemma \ref{le:p-pv}, for
any $y\in \RN$ there exists $\t_y>0$ such that $\tilde w_y=
w_y(\cdot/\t_y)\in \P$. We get
\begin{align*}
|I(\tilde w_y)-b_0| &=|I(\tilde w_y)-I_0(w_y)|
\\
&\le \frac{|\t_y^{N-2}-1|}{2} \irn |\n w|^2 + \frac{\t_y^N}{2}\irn
V(\t_y x+y) w^2
\\
&\quad + |\t_y^N -1| \irn G(w).
\end{align*}
Therefore, by Lemma \ref{le:t}, we infer that
\[
\lim_{|y|\to \infty} I(\tilde w_y)=b_0,
\]
hence $b\le b_0$.
\end{proof}

\begin{lemma}\label{le:pv-p}
Let $z\in \H$ be such that $\irn G(z)>0$. Then there exists
$\bar\t>0$ such that $z^{\bar\t}=z(\cdot/\bar\t)\in \P_0$. In
particular, by ({\bf V6}) this result is true for any $z\in \P$
with $\bar\t \le 1$.
\end{lemma}

\begin{proof}
The first part of the statement follows from the fact that, for
any $z\in \H$ such that $\irn G(z)>0$, certainly there exists
$\bar\theta>0$ such that
    \begin{equation}\label{eq:po}
        \frac{N-2}{2}\irn |\n z|^2 = N\bar\t^2
        \irn G(z).
    \end{equation}
Consider now the case of $z\in \P$. Since
    \begin{equation}\label{eq:po2}
        \frac{N-2}{2}\irn |\n z|^2 + \frac N2 \irn V(x)z^2 +\frac 12 \irn
        (\n V(x)\mid x)z^2 = N \irn G(z),
    \end{equation}
by ({\bf V6}) we have $\irn G(z)>0.$ Let $\bar\theta>0$ such that
\eqref{eq:po} holds. Combining \eqref{eq:po} and \eqref{eq:po2} we
get
\begin{equation}\label{eq:1-t}
\frac 12 \irn \big[ NV(x)+(\n V(x)\mid x) \big] z^2 = N
(1-\bar\t^2)\irn G(z).
\end{equation}
By ({\bf V6}), we get the conclusion.
\end{proof}

Now we can prove Theorem \ref{non}:

\begin{proofnon}
Suppose by contradiction that there exists $z\in \H$ critical
point of the functional $I$ at level $b$: in particular, $z\in \P$
and $I(z)=b$. Let $\t \in (0,1]$ be such that $z^\t\in \P_0$. Let
us show that $\t<1$.
\\
By standard arguments and using the strong maximum principle, we
infer that $z$ does not change sign and so we can assume that
$z>0$. Therefore, by ({\bf V6}) and \eqref{eq:1-t}, we get that
$\t<1$.
\\
By \eqref{eq:iv} and ({\bf V5}), we infer that
\begin{align*}
b&=I(z)= \frac 1N \irn |\n z|^2 -\frac 1{2N} \irn (\n V(x)\mid x)
z^2  \nonumber
\\
&>  \frac {\t^{N-2}}N  \irn |\n z|^2=I_0(z^\t)\ge b_0,
\end{align*}
and we get a contradiction with Lemma \ref{le:<=}.
\end{proofnon}

\begin{proofconon}
If the strict inequality in ({\bf V6}) is satisfied almost
everywhere, then for any $z\in \P$ there exists $\t \in (0,1)$
such that $z^\t\in \P_0$. Arguing as in the proof of Theorem
\ref{non} we conclude.
\end{proofconon}

\begin{remark}
In view of Theorem \ref{non}, the proof of Corollary \ref{co:non} would follow immediately if $\P$ was a natural constraint for the functional $I$. 
\end{remark}


\begin{thebibliography}{99}

\bibitem{BGM}
V. Benci, C.R. Grisanti, A.M. Micheletti,
{\it Existence and non existence of the ground state solution for the nonlinear Schr\"odinger equations with 
$V(\infty)=0$}, Topol. Methods Nonlinear Anal., {\bf 26}, (2005), 203--219.


\bibitem{BL1}
H. Berestycki, P.L. Lions, {\it Nonlinear scalar field equations.
I. Existence of a ground state}, Arch. Rational Mech. Anal., {\bf
82}, (1983), 313--345.


\bibitem{BL2}
H. Berestycki, P.L. Lions, {\it Nonlinear scalar field equations.
II. Existence of infinitely many solutions}, Arch. Rational Mech.
Anal., {\bf 82}, (1983), 347--375.


\bibitem{J}
L. Jeanjean, {\it On the existence of bounded Palais-Smale
sequences and application to a Landesman-Lazer-type problem set on
$\RN$}, Proc. R. Soc. Edinb., Sect. A, Math., {\bf 129},
(1999), 787--809.


\bibitem{J2}
L. Jeanjean, {\it Local condition insuring bifurcation from the
continuous spectrum}, Math. Z., {\bf 232}, (1999), 651--664.


\bibitem{JT}
L. Jeanjean, K. Tanaka, {\it A remark on least energy solutions in
$\RN$}, Proc. Am. Math. Soc., {\bf 131},, (2003) 2399--2408.


\bibitem{JT2}
L. Jeanjean, K. Tanaka, {\it A positive solution for a nonlinear
Schr\"odinger equation on $\RN$}, Indiana Univ. Math. J., {\bf
54}, (2005), 443--464.


\bibitem{R}
P.H. Rabinowitz, \textit{On a class of nonlinear Schr\"odinger
equations}, Z. Angew. Math. Phys., {\bf 43}, (1992), 270--291.


\bibitem{S}
J. Shatah, {\it Unstable ground state of nonlinear Klein-Gordon
equations}, Trans. Amer. Math. Soc., {\bf 290}, (1985), 701--710.


\bibitem{Str}
W.A. Strauss, {\it Existence of solitary waves in higher
dimensions}, Comm. Math. Phys., {\bf 55}, (1977), 149--162.

\end{thebibliography}
\end{document}